\documentclass[12pt,intlim,righttag]{article}
\usepackage{amsmath,amsthm}
\usepackage{amssymb}
\usepackage{amsfonts}

\newcommand{\essinf}{\operatornamewithlimits{ess\,inf}}

\newcommand\beq{\begin{equation}}
\newcommand\eeq{\end{equation}}
\newcommand{\ld}{\lambda}
\newcommand{\vp}{\varphi}

\newcommand{\td}{\tilde}

\newcommand{\lbl}{\label}
\newcommand{\sg}{\sigma}

\newcommand{\la}{\langle}
\newcommand{\ra}{\rangle}

\newcommand{\wt}{\widetilde}
\newcommand{\wh}{\widehat}
\newcommand{\be}{\begin}
\newcommand{\ee}{\end}

\theoremstyle{Theorem}

\theoremstyle{corollary}
\newtheorem{cor}{\ \ \ Corollary}[section]

\theoremstyle{remark}

\theoremstyle{definition}

\numberwithin{equation}{section}

\begin{document}

\title{Mean-variance Hedging Under Partial Information}

\author{M. Mania $^{1),2)}$, R. Tevzadze $^{1),3)}$ and T. Toronjadze $^{1),2)}$}

\date{~}
\maketitle

\begin{center}
$^{1)}$ Georgian American University, Business School,
3, Alleyway II,\\ Chavchavadze Ave. 17, A, Tbilisi, Georgia,

E-mail: toronj333@yahoo.com  \\[2mm]
$^{2)}$ A. Razmadze Mathematical Institute, 1, M. Aleksidze St.,
Tbilisi, Georgia, \\E-mail: mania@rmi.acnet.ge \\[2mm]
$^{3)}$Institute of Cybernetics, 5, S. Euli St., Tbilisi, Georgia,\\
E-mail:tevza@cybernet.ge
\end{center}

\begin{abstract}

We consider the mean-variance hedging problem under partial Information.
The underlying asset price process follows a continuous
semimartingale and strategies have to be constructed when
only part of the information  in the market is available.
We show that the initial  mean variance hedging problem is
equivalent to a new mean variance hedging problem with an
additional correction term, which is formulated in terms
of observable processes. We prove that the value process of the reduced problem
is a square trinomial with coefficients satisfying a triangle
system of backward stochastic differential equations and
the filtered wealth process of the optimal hedging strategy
is characterized as a solution of a linear forward equation.

\bigskip



\end{abstract}

\section{Introduction}

In the problem of derivative pricing and hedging it is usually assumed
that the hedging strategies have to be constructed using all market
information. However, in reality  investors acting in a market have
limited access to the information flow. E.g., an investor may observe
just stock prices,
but stock appreciation rates depend on some unobservable factors;
one may think that stock prices can only be observed at discrete
time instants or with some delay, or an investor would like to
price and hedge a contingent claim whose payoff depends on an unobservable
asset and he observes the prices of an asset correlated with the underlying
asset. Besides, investors may not be able to use  all available
information even if they have access to the full market flow.
In all such cases investors are forced to make decisions based
only on a part of the  market information.

We study a mean-variance hedging problem under partial
information when the asset price process is a continuous semimartingale
and the flow of observable events not necessarily contain
all information on prices of the underlying asset.

We assume that the dynamics of the price process of the asset traded on the
market is described by a continuous semimartingale $S=(S_t,t\in[0,T])$ defined on a filtered probability space
$(\Omega,{\mathcal F}, F=({F}_t,t\in[0,T]),P)$,
satisfying the usual conditions, where ${\cal F} = F_T$
and $T <\infty$ is the fixed time horizon. Suppose that the interest
rate is equal to zero and the asset price process
satisfies the structure condition, i.e.,  the process $S$ admits the decomposition
\begin{equation}\label{str}
S_t=S_0+M_t+\int_0^t\lambda_ud\la M\ra_u,\;\;\;\la\lambda\cdot M\ra_T<\infty\;\;\;a.s.,
\end{equation}
where $M$ is a continuous $F-$local martingale and $\lambda$ is a
$F$-predictable  process.

Let us introduce an additional filtration smaller than F
$$
G_t\subseteq F_t,\;\;\;\;\text{for every}\;\;\; t\in [0,T].
$$
The filtration $G$ represents the information that the hedger has
at his disposal, i.e., hedging strategies have to be constructed using
only information available in $G$.

Let $H$ be a $P$-square integrable $F_T$-measurable
random variable, representing the payoff of a contingent claim
at time $T$.

 We consider the mean-variance hedging problem
\begin{equation}\label{mvh}
\text{to minimize}\;\;\;\;\; E[ (X^{x,\pi}_T-H)^2]\;\;\;\;\text{over
all}\;\;\;\;\pi\in\Pi(G),
\end{equation}
where $\Pi(G)$ is a class of $G$-predictable $S$-integrable
processes. Here
 $X^{x,\pi}_t=x+\int_0^t\pi_udS_u$ is the wealth process starting
from initial capital $x$, determined by the self-financing trading
strategy $\pi\in\Pi(G)$.

In the case $G=F$ of complete information the mean-variance hedging
problem was introduced by F\"ollmer and Sondermann \cite{F-S} in the case
when $S$ is a martingale and then developed by several authors
for price process admitting a trend (see, e.g.,
\cite{D-R}, \cite{H}, \cite{S92},\cite{S94},
\cite{S}, \cite{G-L-Ph}, \cite{H-P-Sc}).

Asset pricing with partial information under various setups has been
considered. The mean-variance hedging problem under partial
information was first studied by Di Masi, Platen and Runggaldier (1995)
when the stock price process is a martingale and the prices
are observed only at discrete time moments. For a general filtrations and
when the asset price process is a martingale this problem was solved
by Schweizer (1994) in terms of $G$-predictable projections.
Pham (2001) considered the mean-variance hedging problem  for a general
semimartingale model, assuming that the
observable filtration contains the augmented filtration $F^S$ generated
by the asset price process $S$
\begin{equation}\label{fsg}
 F^S_t\subseteq G_t ,\;\;\;\;\text{for every}\;\;\; t\in [0,T].
\end{equation}
In this paper, using the variance-optimal martingale
measure with respect to the filtration $G$ and suitable Kunita-Watanabe
decomposition, the theory developed by Gourieroux, Laurent and Pham (1998)
and Rheinl\"ander and Schweizer (1997) to the case of partial information
was extended.

If $F^S_t\subseteq G_t$, the price process is a $G$-semimartingale,
the sharp bracket $\la M\ra$ is $G$-adapted   and the canonical decomposition
of $S$ with respect to the filtration $G$ is of the form
\begin{equation}\label{ino}
S_t=S_0+\int_0^tE(\lambda_u|G_u)d\langle M\rangle_s +\tilde M_t,
\end{equation}
where $\tilde M$ is a $G$-local martingale.

In this case the problem (\ref{mvh}) is equivalent to the problem
\begin{equation}\label{mvhg}
\text{to minimize}\;\;\;\;\; E[ (X^{x,\pi}_T-E(H|G_T))^2]
\;\;\;\;\text{over
all}\;\;\;\;\pi\in\Pi(G)
\end{equation}
which is formulated in $G$-adapted terms, taking in mind
the $G$-decomposition (\ref{ino}) of $S$. Therefore  the problem
(\ref{mvhg})
can be solved as in the case of full information using the dynamic programming
method directly to (\ref{mvhg}) , although one needs
to determine $E(H|G_T)$ and the $G$-decomposition terms of $S$.

If $G$ is not containing $F^S$, then $S$ is not a $G$-semimartingale
and the problem is more involved, although we solve it
under following additional assumptions:

A) $\la M\ra$ and $\lambda$ are $G$-predictable,

B) any $G$- martingale is a $F$-local martingale,

C) the filtration $F$ is continuous, i.e., all $F$- local martingales are continuous,

D) there exists a martingale measure for $S$ that satisfies the Reverse
H\"older condition.

We shall use the notation $\widehat{Y}_t$ for the process
$E(Y_t|G_t)$- the $G$-optional projection of $Y$.
Condition A) implies that
$$
\widehat S_t=E(S_t|G_t)=S_0+\int_0^t\lambda_ud\la M\ra_u+\widehat M_t.
$$

Let
\begin{equation}\label{ht}
H_t=E(H|F_t)=  EH+\int_0^th_udM_u+L_t
\end{equation}
and
\begin{equation}\label{htg}
H_t= EH+\int_0^th^G_udM_u+L^G_t
\end{equation}
be the Galtchouk-Kunita-Watanabe (GKW) decompositions of $H_t=E(H|F_t)$
with respect to local martingales $M$ and $\widehat M$,
where $h, h^G$ are $F$-predictable
process and $L, L^G$ are local martingales strongly
orthogonal to $M$ and $\widehat M$ respectively.

We show (Theorem 3.1) that
the initial mean variance hedging problem (\ref{mvh}) is equivalent
to the problem
to minimize the expression
\begin{equation}\label{opt2i}
E\big[(x+\int_0^T\pi_ud\widehat S_u-
\widehat H_T)^2
+\int_0^T(\pi_u^2(1-\rho^2_u)+2\pi_u\tilde h_u)d\langle M\rangle_u\big],
\end{equation}
over all $\pi\in\Pi(G)$, where
$$
{\widetilde h}_t=
\widehat{h^G_t}\rho^2_t-{\widehat h}_t
\;\;\;\text{and}\;\;\;\;
\rho^2_t=\frac{d\la\widehat M\ra_t}{d\langle M\rangle_t}.
$$

Thus, the problem (\ref{opt2i}), equivalent to (\ref{mvh}), is formulated in terms of $G$-adapted
processes. One can say that (\ref{opt2i}) is the mean variance
hedging problem under complete information with additional correction term and
can be solved as in the case of complete information.

Let us introduce the value process of the problem (\ref{opt2i})

\begin{align}
V^H(t,x)=\underset{\pi\in\Pi(G)}\essinf E\big[(x+\int_t^T\pi_ud\widehat S_u-
\widehat H_T)^2\\
\notag
+\int_t^T[\pi_u^2(1-\rho^2_u)+
2\pi_u\tilde h_u]d\langle M\rangle_u|G_t\big].
\end{align}

We show in section 4 that
the value function of the problem
(\ref{opt2i}) admits a representation
\begin{equation}\label{ut44}
V^H(t,x)=V_t(0)-2V_t(1)x+V_t(2)x^2,
\end{equation}
where the coefficients $V_t(0), V_t(1)$ and $V_t(2)$
 satisfy a triangle
system of backward stochastic differential equations (BSDE).
Besides,
the filtered wealth process of the optimal hedging strategy
is characterized as a solution of the linear forward equation
$$
\widehat X_t^*=x-\int_0^t\frac{\rho_u^2\varphi_u(2)+\lambda_uV_u(2)}
{1-\rho_u^2+\rho_u^2V_u(2)}\widehat X^*_ud\widehat S_u+
$$
\begin{equation}\label{capp}
+\int_0^t\frac{\rho_u^2\varphi_u(1)+\lambda_uV_u(1)+
\tilde h_u}
{1-\rho_u^2+\rho_u^2V_u(2)} d\widehat S_u.
\end{equation}
In the case of complete information ($G=F$) we have $\rho=0$, $\tilde h=0$
and (\ref{capp}) gives equations for the optimal wealth process  and for the coefficients
of value function from \cite{MT7}.

In section 5 we consider a diffusion market model which consists of
two assets $S$ and $Y$, where $S_t$ is a state of a process being controlled
and $Y_t$ is the observation process.  Suppose that $S_t$ and $Y_t$ are
governed by
$$
dS_t=\mu_tdt+\sigma_tdw^0_t,
$$
$$
dY_t=a_tdt+b_tdw_t,
$$
where $w^0$ and $w$ are Brownian motions with correlation $\rho$
and the coefficients
$\mu, \sigma, a$ and $b$ are $F^Y$-adapted. So, in this case $F_t=F_t^{S,Y}$ and
the flow of observable events is $G_t=F^Y_t$. We give in the case of markovian  coefficients
solution of the problem (\ref{mvh}) in terms of parabolic differential equations (PDE) and an
explicit solution when coefficients are constants and the contingent claim is of the form $H={\cal H}(S_T,Y_T)$.

\

\section{Main definitions and auxiliary facts}

\

Denote by ${\mathcal M}^e(F)$ the set of equivalent martingale
measures for $S$, i.e., set of probability measures $Q$
equivalent to $P$ such that $S$ is a  $F$-local martingale under $Q$.

Let
$$
{\mathcal M}^e_2(F)=\{Q\in{\mathcal M}^e(F):
EZ_T^2(Q)<\infty\},
$$
where $Z_t(Q)$ is the density process (with respect to the
filtration $F$) of $Q$ relative to $P$.

{\bf Remark 2.1.}  Since $S$ is continuous, the
 existence of an equivalent martingale measure and the Girsanov
 theorem imply that the structure condition (\ref{str}) is satisfied.

 Note that the density process $Z_t(Q)$
of any element $Q$ of
${\cal M}^e(F)$ is expressed as
an exponential martingale of the form
$$
{\cal E}_t(-\lambda\cdot M+N),
$$
where $N$ is a $F$- local martingale strongly orthogonal to $M$
and  ${\cal E}_t(X)$ is the Doleans-Dade exponential of $X$.

If the local martingale $Z^{min}_t={\cal E}_t(-\lambda\cdot M)$ is
a true martingale, $dQ^{min}/dP= Z^{min}_TdP$ defines
an equivalent probability measure called the minimal martingale
measure for $S$.

Recall that a measure $Q$
satisfies the Reverse H\"older inequality $R_2(P)$ if there exists
a constant $C$ such that
$$
E\big(\frac{Z_T^2(Q)}{Z_\tau^2(Q)}|F_{\tau}\big)\le C, \:\:\:\:\:P-a.s.
$$
for every $F$-stopping time $\tau$.

{\bf Remark 2.2.} If there exists a measure $Q\in{\cal M}^e(F)$ that
satisfies the Reverse H\"older inequality $R_2(P)$, then
according to Kazamaki \cite{Kz} the martingale $M^Q=-\lambda\cdot M+N$ belongs
to the class $BMO$ and hence
$-\lambda\cdot M$ also belongs to $BMO$, i.e.,
\begin{equation}\label{bmo}
E\big(\int_\tau^T\lambda^2_ud\langle M\rangle_u|F_\tau\big)\le const
\end{equation}
for every stopping time $\tau$. Therefore, it follows from Kazamaki \cite{Kz}
that ${\mathcal E}_t(-\lambda\cdot M)$ is
a true martingale. So, condition D) implies that the minimal martingale measure
exists (but $Z^{min}$ is not necessarily square integrable).

For all unexplained notations concerning the martingale theory used below we
refer the reader to \cite{D-M},\cite{L-Sh2},\cite{J}.

Let $\Pi(F)$ be the space of all $F$-predictable $S$-integrable
processes $\pi$ such that the stochastic integral
$$
(\pi\cdot S)_t=\int_0^t\pi_udS_u, t\in[0,T],
$$
is in the ${\cal S}^2$ space of semimartingales , i.e.,
$$
E\big(\int_0^T\pi^2_sd\langle M\rangle_s\big)+
E\big(\int_0^T|\pi_s\ld_s|d\langle M\rangle_s\big)^2<\infty.
$$
Denote by $\Pi(G)$ the subspace of $\Pi(F)$ of $G$-predictable
strategies.

{\bf Remark 2.3.} Since $\lambda\cdot M\in BMO$ (see Remark 2.2),
it follows from the proof of Theorem 2.5 of Kazamaki \cite{Kz}
$$
E\big(\int_0^T|\pi_u\lambda_u|d\langle M\rangle_u\big)^2=
E\la|\pi|\cdot M,|\lambda|\cdot M\ra_T^2
$$
$$
\le 2||\lambda\cdot M||_{\rm BMO} E\int_0^T\pi^2d\langle M\rangle_u<\infty.
$$
Therefore,  under condition D) the strategy $\pi$ belongs to the class $\Pi(G)$
if and only if $E\int_0^T\pi^2_sd\langle M\rangle_s<\infty$.

Define $J_T^2(F)$ and $J_T^2(G)$ as  spaces of terminal values
of stochastic integrals,
i.e.,
$$
J_T^2(F)=\{(\pi\cdot S)_T:\pi\in\Pi(F)\}.
$$
$$
J_T^2(G)=\{(\pi\cdot S)_T:\pi\in\Pi(G)\}.
$$

For convenience we give some assertions from  \cite{D-M-S-S-S} ,
which establishes necessary and sufficient conditions
for the closedness of the space $J_T^2(F)$ in $L^2$.

%{\bf Proposition 2.1.}
\be{prop}\lbl{p2.1}
 Let $S$ be a continuous semimartingale.
Then the following assertions are equivalent:

(1) There is a martingale measure $Q\in{\cal M}^e(F)$ and
$J^2_T(F)$ is closed in $L^2$.

(2) There is a martingale measure $Q\in{\cal M}^e(F)$ that satisfies
the Reverse H\"older condition $R_2(P)$.

(3) There is a constant $C$ such that for all $\pi\in\Pi(F)$ we have
$$
||\sup_{t\le T}(\pi\cdot S)_t||_{L^2(P)}\le C
||(\pi\cdot S)_T||_{L^2(P)}.
$$

(4) There is a constant $c$ such that for every stopping time $\tau$,
every $A\in{\cal F}_\tau$ and for every $\pi\in\Pi(F)$ with
$\pi=\pi I_{]\tau,T]}$ we have
$$
||I_A-(\pi\cdot S)_T||_{L^2(P)}\ge cP(A)^{1/2}.
$$
\ee{prop}

Note that assertion (4) implies that
for every stopping time $\tau$ and for every $\pi\in\Pi(G)$  we have
\begin{equation}\label{vfc}
E\big((1+\int_\tau^T\pi_udS_u)^2/{F}_\tau\big)\ge c.
\end{equation}

 Let us make some remarks on conditions B) and C).

{\bf Remark 2.4.} Conditions $B),C)$ imply that the filtration
$G$ is also continuous. By condition $B$ any $G$-local
martingale is $F$-local martingale,
which are continuous by condition $C)$. Recall that the continuity
of a filtration means that all local martingales with respect to this
filtration are continuous.

{\bf Remark 2.5.} Condition $B)$ is satisfied if and only if
the $\sigma$-algebras $F_t$ and $G_T$
are conditionally independent given $G_t$ for all $t\in[0,T]$
(see Theorem 9.29 from Jacod 1978).

Now we recall some known assertions from the filtering theory.
The following proposition can be proved
similarly to \cite{L-Sh2}.

%{\bf Proposition 2.2.}
\be{prop}\lbl{p2.2}
If conditions $A), B)$ and $C)$ are satisfied,
then for any $F$-local martingale $M$ and any $G$-local martingale $m^G$
\begin{equation}\label{fg}
\widehat M_t=E(M_t|G_t)=
\int_0^tE\big(\frac{d\la M,m^G\ra_u}{d\la m^G\ra_u}|G_u\big)dm^G_u +L^G_t,
\end{equation}
where $L^G$ is a local martingale orthogonal to $m^G$.
\ee{prop}
It follows  from this proposition that for any $G$-predictable, $M$-integrable
process $\pi$ and any $G$-martingale $m^G$
$$
\la\widehat{(\pi\cdot M)},m^G\ra=
\int_0^t\pi_uE\big(\frac{d\la M,m^G\ra_u}{d\la m^G\ra_u}|G_u\big)d\la m^G\ra_u=
$$
$$
=\int_0^t\pi_ud\la\widehat M,m^G\ra_u=\la\pi\cdot\widehat M,m^G\ra_t.
$$
Hence, for any $G$-predictable, $M$-integrable
process $\pi$
\begin{equation}\label{mg}
\widehat{(\pi\cdot M)_t}=E\big(\int_0^t\pi_sdM_s|G_t)=\int_0^t\pi_sd\widehat M_s.
\end{equation}
Since $\pi, \lambda$ and $\langle M\rangle$ are $G$-predictable, from (\ref{mg}) we have
 \begin{equation}
\widehat{(\pi\cdot S)_t}=E\big(\int_0^t\pi_udS_u|G_t)=
\int_0^t\pi_ud\widehat S_u,
\end{equation}
where
$$
 \widehat S_t=S_0+\int_0^t\lambda_ud\langle M\rangle_u+\widehat M_t.
$$

\section{Separation principle. The optimality principle}

Let us introduce the value function of the problem (\ref{mvh}) defined as
\begin{equation}
U^H(t,x)=\underset{\pi\in\Pi(G)}{\essinf}
E\big((x+\int_t^T\pi_udS_u-H )^2|{G}_t\big).
\end{equation}
By GKW decomposition
\begin{equation}\label{hhh2}
H_t=E(H|F_t)=  EH+\int_0^th_udM_u+L_t
\end{equation}
for a $F$-predictable, $M$-integrable
process $h$ and a local martingale $L$ strongly
orthogonal to $M$.
We shall use also the GKW decompositions of $H_t=E(H|F_t)$
with respect to the local martingale $\widehat M$
\begin{equation}\label{htg}
H_t= EH+\int_0^th^G_ud\wh{M}_u+L^G_t
\end{equation}
where $h^G$ is a $F$-predictable
process and $L^G$ is a $F$- local martingale strongly
orthogonal to $\widehat M$.

It follows from  Proposition \ref{p2.2} ( applied for $m^G=\widehat M$)
and  Lemma \ref{l2.1} that

\begin{equation}\label{hmbar}
\la E(H|G_.),\widehat M\ra_t=\int_0^tE(h^G_u|G_u)d\la\widehat M\ra_u=
\int_0^t\widehat{h_u^G}\rho^2_ud\la M\ra_u.
\end{equation}

We shall use the  notation
\begin{equation}\label{ro}
\widetilde h_t=
\widehat{h_t^G}\rho^2_t-\widehat h_t.
\end{equation}

Note that $\widetilde h$ belongs to the class $\Pi(G)$ by Lemma \ref{l3.1}.

Let us introduce now a new optimization problem,
equivalent to the initial mean variance hedging problem (\ref{mvh}),
to minimize the expression
\begin{equation}\label{opt2}
E\big[(x+\int_0^T\pi_ud\widehat S_u-
\widehat H_T)^2
+\int_0^T(\pi_u^2(1-\rho^2_u)+2\pi_u\tilde h_u)d\langle M\rangle_u\big],
\end{equation}
over all $\pi\in\Pi(G)$. Recall that $\widehat S_t= E(S_t|G_t)= S_0+
\int_0^t\lambda_ud\langle M\rangle_u+\widehat M_t$.

%Theorem 3.1.
\be{Th}\lbl{t3.1}
Let conditions $A), B)$ and $C)$ be
satisfied. Then the initial mean-variance hedging problem
(\ref{mvh}) is equivalent to the problem (\ref{opt2}). In
particular, for any $\pi\in\Pi(G)$ and $t\in [0,T]$
\begin{equation}\label{eqv}
E\big[(x+\int_t^T\pi_udS_u-H)^2|G_t\big]=
E\big[(H-\widehat H_T)^2|G_t\big]
\end{equation}
$$
+E\big[(x+\int_t^T\pi_ud\widehat S_u-\widehat H_T)^2
+\int_t^T(\pi_u^2(1-\rho^2_u)+2\pi_u\tilde h_u)d\langle M\rangle_u|G_t\big].
$$
\ee{Th}
{\it Proof.} We have
$$
E\big[(x+\int_t^T\pi_udS_u-H)^2|G_t\big]=
E\big[\big(x+\int_t^T\pi_ud\widehat S_u-H+
\int_t^T\pi_ud(M_u-\widehat M_u)\big)^2|G_t\big]
$$
$$
=E\big[\big(x+\int_t^T\pi_ud\widehat S_u-H\big)^2|G_t\big]+
 2E\big[\big(x+\int_t^T\pi_ud\widehat S_u-H\big)
\big( \int_t^T\pi_ud(M_u-\widehat M_u)\big)|G_t\big]
$$
\begin{equation}\label{i}
+E\big[\big(\int_t^T\pi_ud(M_u-\widehat M_u)\big)^2|G_t\big]=
I_1+2I_2+I_3.
\end{equation}
It is evedent that
\begin{equation}\label{i1}
I_1=E\big[(x+\int_t^T\pi_ud\widehat S_u-\widehat H_T)^2|G_t\big]+
E\big[(H-\widehat H_T)^2|G_t\big].
\end{equation}
Since $\pi,\lambda$ and $\langle M\rangle$ are $G_T$-measurable and the $\sigma$-algebras $F_t$
and $G_T$
are conditionally independent given $G_t$ (see Remark 2.5), it follows from
equation (\ref{mg}) that
$$
E\big[\int_t^T\pi_u\lambda_ud\langle M\rangle_u \int_t^T\pi_ud(M_u-\widehat M_u)|G_t\big]
=E\big[\int_t^T\pi_u\lambda_ud\langle M\rangle_u \int_0^T\pi_ud(M_u-\widehat M_u)|G_t\big]
$$
$$
-E\big[\int_t^T\pi_u\lambda_ud\langle M\rangle_u \int_0^t\pi_ud(M_u-\widehat M_u)|G_t\big]
=E\big[\int_t^T\pi_u\lambda_ud\langle M\rangle_u E(\int_0^T\pi_ud(M_u-\widehat M_u)|G_T)|G_t\big]
$$
\begin{equation}\label{i2}
-E\big[\int_t^T\pi_u\lambda_ud\langle M\rangle_u |G_t\big]E\big[\int_0^t\pi_ud(M_u-\widehat M_u)|G_t\big]=
0
\end{equation}

On the other hand using decomposition (\ref{hhh2}), equality
(\ref{hmbar}), properties of
square characteristics of martingales and the projection
theorem we obtain
$$
E\big[H\int_t^T\pi_ud(M_u-\widehat M_u)|G_t\big]=
E\big[H\int_t^T\pi_udM_u|G_t\big]-E\big[\widehat H_T\int_t^T\pi_ud\widehat M_u|G_t\big]
$$
$$
=E\big[\int_t^T\pi_ud\la M,E(H|F_\cdot)\ra_u|G_t\big]-E\big[\int_t^T\pi_ud\la\widehat H,
\widehat M\ra_u|G_t\big]
$$
$$
=E\big[\int_t^T\pi_uh_ud\langle M\rangle_u|G_t]-E\big[\int_t^T\pi_u
\widehat{h_u^G}\rho^2_ud\langle M\rangle_u|G_t\big]=
$$
\begin{equation}\label{i2h}
E\big[\int_t^T\pi_u(\widehat h_u-\widehat{h_u^G}\rho^2_u)d\la M\ra_u|G_t\big]
= -E\big[\int_t^T\pi_u\widetilde h_ud\langle M\rangle_u|G_t\big].
\end{equation}

Finally, it is easy to verify that
$$
2E\big[\int_t^T\pi_u\widehat M_u\int_t^T\pi_ud(M_u-\widehat M_u)|G_t\big]+
E\big[\big(\int_t^T\pi_ud(M_u-\widehat M_u)\big)^2|G_t\big]=
$$
$$
E\big[(\int_t^T\pi_u^2d\langle M\rangle_u-\int_t^T\pi_u^2d\la\widehat M\ra_u)|G_t\big]=
$$
\begin{equation}\label{i3}
=E\big[\int_t^T\pi_u^2(1-\rho^2_u)d\langle M\rangle_u|G_t\big].
\end{equation}
Therefore equations (\ref{i}), (\ref{i1}),(\ref{i2}), (\ref{i2h}),
and (\ref{i3}) imply the validity of equality (\ref{eqv}). \qed

Thus, it follows from Theorem \ref{t3.1} that the optimization problems
(\ref{mvh}) and (\ref{opt2}) are equivalent. Therefore it is sufficient to
solve the
problem (\ref{opt2}), which is formulated in terms of $G$-adapted
processes. One can say that (\ref{opt2}) is a mean variance
hedging problem under complete information with correction term and
can be solved as in the case of complete information.

Let us introduce the value process of the problem (\ref{opt2})

$$
V^H(t,x)= \essinf_{\pi\in\Pi(G)} E\big[(x+\int_t^T\pi_ud\widehat S_u-
\widehat H_T)^2+
$$
\begin{equation}\label{vth}
+\int_t^T[\pi_u^2(1-\rho^2_u)+
2\pi_u\tilde h_u]d\langle M\rangle_u|G_t\big].
\end{equation}

It follows from Theorem \ref{t3.1} that
\begin{equation}\label{uv}
U^H(t,x)=V^H(t,x)+
E\big[(H-\widehat H_T)^2|G_t\big].
\end{equation}

The optimality principle takes in this case the following form

%{Proposition 3.1}
\be{prop}\lbl{p3.1}
(Optimality principle). Let conditions $A, B)$ and $C)$
 be satisfied. Then

a) For all $x\in R$, $\pi\in\Pi(G)$ and $s\in[0,T]$ the process
$$
V^H(t,x+\int_s^t\pi_ud\widehat S_u)+  \int_s^t[\pi_u^2(1-\rho^2_u)+
2\pi_u\tilde h_u)]d\langle M\rangle_u
$$
is a submartingale on $[s,T]$,
admitting an RCLL modification.

b) $\pi^*$ is optimal if and only if the process

$$
V^H(t,x+\int_s^t\pi^*_ud\widehat S_u)+  \int_s^t[(\pi^*_u)^2(1-\rho^2_u)+
2\pi^*_u\tilde h_u]d\langle M\rangle_u
$$
is a martingale.
\ee{prop}
This assertion can be proved in a standard manner (see, e.g., \cite{El-Q}, \cite{Kr}).
The proof more adapted to this case one can see in \cite{MT7}.

Let
$$
V(t, x)=\underset{\pi\in\Pi(G)}{\essinf}
E\big[(x+\int_t^T\pi_ud\widehat S_u)^2+
\int_t^T\pi_u^2(1-\rho^2_u)d\langle M\rangle_u|G_t\big].
$$
and
$$
V_t(2)=\underset{\pi\in\Pi(G)}{\essinf}
E\big[(1+\int_t^T\pi_ud\widehat S_u)^2+
\int_t^T\pi_u^2(1-\rho^2_u)d\langle M\rangle_u|G_t\big].
$$
It is evident that $V(t,x)$ (resp. $V_2(t)$) is the
value process of the optimization problem (\ref{opt2})
in the case $H=0$ (resp. $H=0$ and $x=1$), i.e.,
$$
V(t,x)=V^0(t,x)\;\;\;\text{and}\;\;\; V_2(t)=V^0(t,1).
$$

Since $\Pi(G)$
is a cone, we have that
$$
V(t,x)=
x^2\underset{\pi\in\Pi(G)}{\essinf}
E\big[(1+\int_t^T\frac{\pi_u}{x}d\widehat S_u)^2+
$$
\begin{equation}\label{vx2}
+\int_t^T\big(\frac{\pi_u}{x}\big)^2(1-\rho^2_u)
d\langle M\rangle_u|G_t\big]= x^2 V_2(t).
\end{equation}

Therefore from Proposition \ref{p3.1} and equality (\ref{vx2})
we have the following

%{\bf Corollary 3.1.}
\be{cor}\lbl{c3.1}
 a) The process
$$
V_2(t)(1+\int_s^t\pi_ud\widehat S_u)^2+
\int_s^t(\pi_u)^2(1-\rho^2_u)d\langle M\rangle_u,
$$
$t\ge s)$ is a submartingale
for all $\pi\in\Pi(G)$ and $s\in[0,T]$.

b) $\pi^*$ is optimal
iff
$$
V_2(t)(1+\int_s^t\pi^*_ud\widehat S_u)^2+
\int_s^t(\pi^*_u)^2(1-\rho^2_u)d\langle M\rangle_u,
$$
$ t\ge s,$ is a martingale.
\ee{cor}
Note that in the case $H=0$ from  Theorem 3.1
we have

\begin{equation}\label{eqv0}
E\big[(1+\int_t^T\pi_udS_u)^2|G_t\big)=
\end{equation}
$$
E\big[(1+\int_t^T\pi_ud\widehat S_u)^2
+\int_t^T\pi_u^2(1-\rho^2_u)d\langle M\rangle_u|G_t\big]
$$
and, hence
\begin{equation}\label{vequ}
V_2(t)=U^0(t,1).
\end{equation}
%{\bf Lemma 3.2.}
\be{lem}\lbl{l3.2}
Let conditions $A)-D)$  be satisfied.
Then there is a constant $1\ge c>0$ such that
$V_t(2)\ge c$ for all $t\in[0,T]$ a.s. and
\begin{equation}\label{rgec}
1-\rho^2_t+\rho^2_tV_t(2)\ge c\;\;\;\;\;\mu^{\langle M\rangle} a.e.
\end{equation}
\ee{lem}
{\it Proof.}  Let
$$
V_t^F(2)=\underset{\pi\in\Pi(F)}{\essinf}
E\big[(1+\int_t^T\pi_udS_u)^2|F_t\big].
$$
It follows from assertion 4) of Proposition \ref{p2.1} that
there is a constant $c>0$ such that
$V^F_t(2)\ge c$ for all $t\in[0,T]$ a.s.. Note that $c\le 1$ since
$V^F\le 1$. Then by (\ref{vequ})
$$
V_t(2)=U^0(t,1)=\underset{\pi\in\Pi(G)}{\essinf}
E\big[(1+\int_t^T\pi_udS_u)^2|G_t\big]=
$$
$$
V_t(2)=\underset{\pi\in\Pi(G)}{\essinf}
E\big[E((1+\int_t^T\pi_udS_u)^2|F_t)|G_t\big]\ge
$$
$$
\ge V_t^F(2)\ge c.
$$
Therefore, since $\rho^2_t\le 1$ by Lemma \ref{l2.1},
$$
1-\rho^2_t+\rho^2_tV_t(2)\ge 1-\rho^2_t+\rho^2_tc\ge
\inf_{r\in[0,1]}(1-r+rc)=c.
$$

\

\section{BSDEs for the value process}

\

Let us consider the semimartingale backward equation
\begin{equation}\label{sbe}
Y_t=Y_0+\int_0^tf(u, Y_u, \psi_u)d\la m\ra_u+\int_0^t\psi_udm_u+L_t
\end{equation}
with the boundary condition
\begin{equation}\label{bc}
Y_T=\eta,
\end{equation}
where $\eta$ is an integrable $G_T$-measurable random variable,
$f:\Omega\times [0,T]\times R^2\to R$ is ${\cal P}\times{\cal B}(R^2)$
measurable and $m$ is a local martingale.
A solution of (\ref{sbe})-(\ref{bc}) is a triple $(Y,\psi,L)$, where
$Y$ is a special semimartingale, $\psi$ is a predictable $m$-integrable
process and $L$ a local martingale strongly orthogonal to $m$. Sometimes
we call $Y$ alone the solution of (\ref{sbe})-(\ref{bc}), keeping in mind
that $\psi\cdot m+L$ is the martingale part of $Y$.

Backward stochastic differential equations  have been introduced in
\cite{B} for the linear case as the equations for the adjoint process in the stochastic
maximum principle.The semimartingale backward equation, as a stochastic version of the Bellman
equation in an optimal control problem, was first derived in \cite{Ch}. The BSDEs with more
general nonlinear generators  was introduced in \cite{Par-P} for the case of Brownian filtration,
where an existence and uniqueness of a solution of BSDEs with generators satisfying the global
Lifschitz condition was established. These results were generalized for
generators with quadratic growth in  \cite{Kob}, \cite{LS} for BSDEs driven by a Brownian motion and
in \cite{Mrl}, \cite{Tev} for BSDEs driven by martingales. But conditions imposed in these papers
are too restrictive for our needs. We prove here existence and uniqueness of a solution by directly
showing that the unique solution of the BSDE we consider is the value of the problem.

In this section we  characterize optimal strategies in terms of solutions of suitable Semimartingale Backward Equations.

%Theorem 4.1.
\be{Th}\lbl{t4.1}
Let $H$ be a square integrable $F_T$-measurable
random variable and let conditions $A), B), C)$ and $D)$ be satisfied.
 Then the value function of the problem
(\ref{opt2}) admits a representation
\begin{equation}\label{ut44}
V^H(t,x)=V_t(0)-2V_t(1)x+V_t(2)x^2,
\end{equation}
where the processes $V_t(0), V_t(1)$ and $V_t(2)$ satisfy the
following system of backward equations
$$
Y_t(2)= Y_0(2)+\int_0^t \frac{\big(\psi_s(2)\rho_s^2+\lambda_sY_s(2)\big)^2}
{1-\rho_s^2+\rho_s^2Y_s(2)}d\langle M\rangle_s
$$
\begin{equation}\label{v2}
+\int_0^t\psi_s(2)d\widehat M_s+ L_t(2)\;\;\;\;Y_T(2)=1,
\end{equation}
$$
Y_t(1)= Y_0(1)+\int_0^t \frac{\big(\psi_s(2)\rho_s^2+\lambda_sY_s(2)\big)
\big(\psi_s(1)\rho_s^2+\lambda_sY_s(1)-\tilde h_s\big)}
{1-\rho_s^2+\rho_s^2Y_s(2)}d\langle M\rangle_s
$$
\begin{equation}\label{v1}
+\int_0^t\psi_s(1)d\widehat M_s+ L_t(1),\;\;\;\;Y_T(1)=E(H|G_T),
\end{equation}
$$
Y_t(0)= Y_0(0)+\int_0^t \frac{\big(\psi_s(1)\rho_s^2+\lambda_sY_s(1)-\tilde h_s\big)^2}
{1-\rho_s^2+\rho_s^2Y_s(2)}d\langle M\rangle_s
$$
\begin{equation}\label{v0}
+\int_0^t\psi_s(0)d\widehat M_s+ L_t(0),\;\;Y_T(0)=E^2(H|G_T),
\end{equation}
where $L(2), L(1)$ and $L(0)$ are $G$-local martingales
orthogonal to $\widehat M$.

Besides the  optimal  filtered wealth process
$\widehat X_t^{x,\pi^*}=x+\int_0^t\pi^*_ud\widehat S_u$
 is a solution of the
linear equation
$$
\widehat X_t^*=x-\int_0^t\frac{\rho_u^2\varphi_u(2)+\lambda_uV_u(2)}
{1-\rho_u^2+\rho_u^2V_u(2)}\widehat X^*_ud\widehat S_u+
$$
\begin{equation}\label{cap}
+\int_0^t \frac{\varphi_u(1)\rho_u^2+\lambda_uV_u(1)-\tilde h_u}
{1-\rho_u^2+\rho_u^2V_u(2)}
 d\widehat S_u.
\end{equation}
\ee{Th}
{\it Proof.} Similarly to the case of complete information one can show
that the optimal strategy exists and that $V^H(t,x)$ is a square
trinomial of the form (\ref{ut44}) (see, e.g.,  \cite{MT7}).
More precisely the space  of stochastic integrals
$$
J^2_T(G)=\{(\pi\cdot S)_T:\pi\in\Pi(G)\}
$$
is closed by Proposition \ref{p2.1} and condition A). Hence there exists optimal
strategy $\pi^*(t,x)\in \Pi(G)$ and
$U^H(t,x)=E[|H-x-\int_t^T\pi_u^*(t,x)dS_u|^2|{\mathcal F}_t].$ Since
$\int_t^T\pi_u^*(t,x)dS_u$ coincides with the orthogonal
projection of $H-x\in L^2$ on the closed subspace of stochastic
integrals, then the optimal strategy is linear with respect to
$x$, i.e., $\pi_u^*(t,x)=\pi_u^0(t)+x{\pi_u^1}(t).$ This implies
that the value function $U^H(t,x)$  is a square trinomial. It follows
from the equality
(\ref{uv})that $V^H(t,x)$ is also a square trinomial
and it admits the representation (\ref{ut44}).

Let us show that $V_t(0), V_t(1)$ and $V_t(2)$ satisfy the
system (\ref{v2})-(\ref{v0}). It is evident that
$$
V_t(0)=V^H(t,0)=\underset{{\pi\in\Pi(G)}}
\essinf E\big[(\int_t^T\pi_ud\widehat S_u-
\widehat H_T)^2
$$
\begin{equation}\label{v0t}
+\int_t^T[\pi_u^2(1-\rho^2_u)+
2\pi_u\tilde h_u]d\langle M\rangle_u|G_t\big]
\end{equation}
and
$$
V_t(2)=V^0(t,1)=\underset{{\pi\in\Pi(G)}}\essinf E\big[(1+\int_t^T\pi_ud\widehat S_u)^2
$$
\begin{equation}\label{v2t}
+\int_t^T\pi_u^2(1-\rho^2_u)d\langle M\rangle_u|G_t\big].
\end{equation}
Therefore, it follows from the optimality principle (taking $\pi=0$)
that  $V_t(0)$ and $V_t(2)$ are RCLL $G$-submaringales  and
$$
V_t(2)\le E(V_2(T)|G_t)\le1,
$$
$$
V_0(t)\le E(E^2(H|G_T)|G_t)\le E(H^2|G_t).
$$
Since
\begin{equation}\label{v1t}
V_t(1)=\frac{1}{2}(V_t(0)+V_t(2)-V^H(t,1)),
\end{equation}
the process $V_t(1)$
is also a special semimartingale and since
$V_t(0)-2V_t(1)x+V_t(2)x^2=V^H(t,x)\ge0$ for all $x\in R$, we have that
$V_t^2(1)\le V_t(0)V_t(2)$, hence
$$
V_t^2(1)\le E(H^2|G_t).
$$

Expressions (\ref{v0t}), (\ref{v2t}) and (\ref{vth})  imply that
$V_T(0)=E^2(H|G_T)$, $V_T(2)=1$ and $V^H(T,x)=(x-E(H|G_T))^2$.
Therefore from (\ref{v1t}) we have $V_T(1)=E(H|G_T)$ and $V(0), V(1), V(2)$
satisfy the boundary conditions.

Thus, the coefficients $V_t(i), i=0,1,2$ are special
semimartingales and they admit the decomposition
\begin{equation}\label{gkw2}
V_t(i)= V_0(i) + A_t(i) + \int_0^t\varphi_s(i)d\widehat M_s+ m_t(i),\;\;\;
i=0,1,2,
\end{equation}
where $m(0), m(1), m(2)$ are $G$-local martingales strongly orthogonal to
$\widehat M$.

There exists an increasing continuous $G$-predictable process $K$
such that
$$
\langle M\rangle_t=\int_0^t\nu_udK_u,\;\;\;\;A_t(i)=\int_0^ta_u(i)dK_u,\;\;\;
i=0, 1, 2,
$$
where $\nu$ and $a(i),i=0,1,2,$ are $G$-predictable processes.

Let $\widehat X_{s,t}^{x,\pi}\equiv x+\int_s^t\pi_ud\widehat S_u$
and
$$
Y^{x,\pi}_{s,t}\equiv V^H(t, \widehat X_{s,t}^{x,\pi})+
\int_s^t[\pi_u^2(1-\rho^2_u)+
2\pi_u\tilde h_u)]d\langle M\rangle_u.
$$
Then using (\ref{ut44}), (\ref{gkw2}) and the It\^o formula
for any $t\ge s$ we have
$$
(\widehat X_{s,t}^{x,\pi})^2=x + \int_s^t[2\pi_u\lambda_u
\widehat X_{s,u}^{x,\pi}+\pi_u^2\rho^2_u]d\langle M\rangle_u+
$$
\begin{equation}\label{x2}
+2\int_s^t\pi_u\widehat X_{s,u}^{x,\pi}d\widehat M_u
\end{equation}
and
$$
Y^{x,\pi}_{s,t}- V^H(s,x)=
\int_s^t[(\widehat X_{s,u}^{x,\pi})^2a_u(2)-
2\widehat X_{s,u}^{x,\pi}a_u(1)
+a_u(0)]dK_u+
$$
$$
\int_s^t\big[\pi_u^2(1-\rho^2_u+\rho^2_uV_{u-}(2))+
2\pi_u\widehat X_{s,u}^{x,\pi}(\lambda_uV_{u-}(2)
+\varphi_u(2)\rho^2_u)-
$$
\begin{equation}\label{ito}
-2\pi_u(V_{u-}(1)\lambda_u
+\varphi_u(1)\rho^2_u-\tilde h_u)]\nu_udK_u+m_t-m_s,
\end{equation}
where $m$ is a local martingale.

Let
$$
G(\pi,x)=G(\omega,t,\pi,x)=
\pi^2(1-\rho^2_u+\rho^2_uV_{u-}(2))+
2\pi x(\lambda_uV_{u-}(2)
+\varphi_u(2)\rho^2_u)-
$$
$$
-2\pi(V_{u-}(1)\lambda_u
+\varphi_u(1)\rho^2_u-\tilde h_u)
$$

It follows from the optimality principle
that for each $\pi\in\Pi(G)$ the process
$$
\int_s^t[(\widehat X_{s,u}^{x,\pi})^2a_u(2)-
2\widehat X_{s,u}^{x,\pi}a_u(1))+a_u(0)]dK_u+
$$
\begin{equation}\label{optge}
+\int_s^tG(\pi_u,\widehat X_{s,u}^{x,\pi})\nu_udK_u
\end{equation}
is increasing for any $s$ on $s\le t\le T$ and for the optimal strategy
$\pi^*$ we have the equality
$$
  \int_s^t[(\widehat X_{s,u}^{x,\pi^*})^2a_u(2)-
2\widehat X_{s,u}^{x,\pi^*}a_u(1)+a_u(0)]dK_u=
$$
\begin{equation}\label{opteq}
-\int_s^tG(\pi^*_u,\widehat X_{s,u}^{x,\pi^*})\nu_udK_u.
\end{equation}

Since $\nu_udK_u=d\langle M\rangle_u$ is continuous, without loss of generality
one can assume that the process $K$ is continuous (see \cite{MT7} for details).
Therefore, taking in (\ref{optge})
$\tau_s(\varepsilon)
=\inf\{t\ge s:K_t-K_s\ge\varepsilon\}$
instead of $t$ we have that for any $\varepsilon>0$ and $s\ge0$
$$
\frac{1}{\varepsilon}\int_s^{\tau_s(\varepsilon)}
 [(\widehat X_{s,u}^{x,\pi})^2a_u(2)-
2\widehat X_{s,u}^{x,\pi}a_u(1)+a_u(0)]dK_u\ge
$$
\begin{equation}\label{limit}
-\frac{1}{\varepsilon}\int_s^{\tau_s(\varepsilon)}
G(\pi_u,\widehat X_{s,u}^{x,\pi})\nu(u)dK_u.
\end{equation}
Passing to the limit in (\ref{limit}) as $\varepsilon\to 0$, from Proposition B
of \cite{MT7} we obtain that
$$
x^2a_u(2)-
2xa_u(1)+a_u(0)
\ge -G(\pi_u, x)\nu_u\;\;\;\;\;\;\mu^{K}-a.e.
$$
for all $\pi\in\Pi(G).$   Similarly from (\ref{opteq}) we have
that $\mu^{K}$-a.e.
$$
x^2a_u(2)-
2xa_u(1)+a_u(0)
= -G(\pi_u, x)\nu_u\;\;\;\;\;\;
$$
and hence
\begin{equation}\label{ess}
x^2a_u(2)-
2xa_u(1)+a_u(0)= -\nu_u\underset{\pi\in\Pi(G)}{\essinf}
G(\pi_u, x).
\end{equation}
The infinum in (\ref{ess}) is attained
for the strategy
\begin{equation}\label{phat}
\hat\pi_t= \frac{V_t(1)\lambda_t+\varphi_t(1)\rho_t^2-\tilde h_t-
x(V_t(2)\lambda_t+\varphi_t(2)\rho_t^2)}
{1-\rho^2_t+\rho_t^2V_t(2)}.
\end{equation}
From here we can conclude that
$$
\underset{\pi\in\Pi(G)}{\essinf}
G(\pi_t, x)\ge G(\hat\pi_t,x)=
$$
\begin{equation}\label{infi}
=-\frac{(V_t(1)\lambda_t+\varphi_t(1)\rho_t^2-\tilde h_t-
x(V_t(2)\lambda_t+\varphi_t(2)\rho_t^2))^2}
{1-\rho^2_t+\rho_t^2V_t(2)}.
\end{equation}
Let $\pi^n_t=I_{[0,\tau_n[}(t)\hat\pi_t$, where
$\tau_n=\inf\{t:|V_t(1)|\ge n\}$.

It follows from Lemma \ref{l3.1}, Lemma \ref{l3.2} and Lemma \ref{la4} that
$\pi^n\in\Pi(G)$ for every $n\ge1$ and hence
$$
\underset{\pi\in\Pi(G)}{\essinf}
G(\pi_t, x)\le G(\pi^n_t,x)
$$
for all $n\ge1$. Therefore
\begin{equation}\label{hat}
\underset{\pi\in\Pi(G)}{\essinf}
G(\pi_t, x)\le\lim_{n\to\infty}G(\pi^n_t,x)=G(\hat\pi_t,x).
\end{equation}
Thus (\ref{ess}), (\ref{infi}) and (\ref{hat}) imply that
$$
x^2a_t(2)-
2xa_t(1)+a_t(0)=
$$
\begin{equation}\label{trin}
=\nu_t\frac{(V_t(1)\lambda_t+\varphi_t(1)\rho_t^2-\tilde h_t-
x(V_t(2)\lambda_t+\varphi_t(2)\rho_t^2)^2}
{1-\rho^2_t+\rho_t^2V_t(2)}, \;\;\;\;\mu^{K}\;\;\;a.e.
\end{equation}
and equalizing the coefficients of square
trinomials in (\ref{trin}) (and integrating with respect to
$dK$) we obtain that
\begin{equation}\label{at2}
A_t(2)= \int_0^t\frac{\big(\varphi_s(2)\rho_s^2+
\lambda_sV_s(2)\big)^2}
{1-\rho_s^2+\rho_s^2V_s(2)}d\langle M\rangle_s,
\end{equation}
\begin{equation}\label{at1}
A_t(1)=\int_0^t \frac{\big(\varphi_s(2)\rho_s^2+\lambda_sV_s(2)\big)
\big(\varphi_s(1)\rho_s^2+\lambda_sV_s(1)-\tilde h_s\big)}
{1-\rho_s^2+\rho_s^2V_s(2)}d\langle M\rangle_s,
\end{equation}
\begin{equation}\label{at0}
A_t(0)=
\int_0^t \frac{
\big(\varphi_s(1)\rho_s^2+\lambda_sV_s(1)-\tilde h_s\big)^2}
{1-\rho_s^2+\rho_s^2V_s(2)}d\langle M\rangle_s,
\end{equation}
which, together with (\ref{gkw2}), implies that the triples $(V(i),\varphi(i),
m(i))$, $ i=0, 1, 2,$
satisfy the system (\ref{v2})-(\ref{v0}).

Note that $A(0)$ and $A(2)$ are integrable increasing processes and relations
(\ref{at2}) and (\ref{at0}) imply that the strategy $\hat\pi$ defined by
(\ref{phat}) belongs to the class $\Pi(G)$.

Let us show now that if the  strategy
$\pi^*\in\Pi(G)$ is optimal then the corresponding filtered wealth
process $\widehat X_t^{\pi^*}=x+\int_0^t\pi_u^*d\widehat S_u$ is
a solution of equation (\ref{cap}).

By the optimality principle the process
$$
Y^{\pi^*}_t=V^H(t,\widehat X^{\pi^*}_t)+  \int_0^t[(\pi^*_u)^2(1-\rho^2_u)+
2\pi^*_u\tilde h_u]d\langle M\rangle_u
$$
is a martingale. Using the It\^o formula we have
$$
 Y_t^{\pi^*}= \int_0^t(\widehat X_{u}^{\pi^*})^2dA_u(2)-
2\int_0^t\widehat X_{u}^{\pi^*}dA_u(1)+A_t(0)
$$
\begin{equation}\label{optpr}
+\int_0^tG(\pi^*_u,\widehat X_{u}^{\pi^*})d\langle M\rangle_u+ N_t,
\end{equation}
where $N$ is a martingale. Therefore applying equalities (\ref{at2}),(\ref{at1}) and (\ref{at0})
we obtain that
$$
 Y_t^{\pi^*}= \int_0^t\big(\pi^*_u - \frac{V_u(1)\lambda_u+\varphi_u(1)\rho_u^2-\tilde h_u}
{1-\rho^2_u+\rho_u^2V_u(2)}
$$
\begin{equation}\label{optpr}
+\widehat X_u^{\pi^*}\frac{V_u(2)\lambda_u+\varphi_u(2)\rho_u^2}
{1-\rho^2_u+\rho_u^2V_u(2)}\big)^2 (1-\rho^2_u+\rho_u^2V_u(2))d\langle M\rangle_u+ N_t,
\end{equation}
which implies that $\mu^{\langle M\rangle}$-a.e.
$$
\pi^*_u = \frac{V_u(1)\lambda_u+\varphi_u(1)\rho_u^2-\tilde h_u}
{1-\rho^2_u+\rho_u^2V_u(2)}-
\widehat X_u^{\pi^*}\frac{(V_u(2)\lambda_u+\varphi_u(2)\rho_u^2)}
{1-\rho^2_u+\rho_u^2V_u(2)}
$$
Integrating both parts of this equality with respect to $d\widehat S$ (and adding then $x$ to
the both parts) we obtain that $\widehat X^{\pi^*}$ satisfies equation (\ref{cap}).\qed

The uniqueness of the system  (\ref{v2})-(\ref{v0}) we shall
prove under following condition $D^*)$, stronger than condition D).

Assume that

$D^*)$
$$
\int_0^T\frac{\lambda^2_u}{\rho^2_u}d\langle M\rangle_u\le C.
$$
Since $\rho^2\le1$ ( Lemma \ref{l2.1}), it follows from $D^*)$ that the
mean-variance tradeoff of $S$ is bounded, i.e.,
$$
 \int_0^T\lambda^2_ud\langle M\rangle_u\le C,
$$
which implies that the minimal martingale measure for $S$ exists and
satisfies the Reverse-H\"older condition $R_2(P)$. So, condition $D^*)$
implies condition D). Besides it follows from the condition $D^*)$ that
the minimal martingale measure ${\widehat Q}^{min}$ for $\widehat S$
$$
d{\widehat Q}^{min}={\cal E}_T(-\frac{\lambda}{\rho^2}\cdot\widehat M)
$$
also exists and satisfies  Reverse-H\"older condition.

Recall that the process $Z$ belongs to the class $D$ if the family of
random variables $Z_\tau I_{(\tau\le T)}$ for all stopping times
$\tau$ is uniformly integrable.

%{\bf Theorem 4.2}
\be{Th}\lbl{t4.2}
Let conditions $A), B), C)$ and $D^*)$ be satisfied.
If a triple $(Y(0),Y(1),Y(2))$, where $Y(0)\in D, Y^2(1)\in D$ and
$c\le Y(2)\le C$ for some constants $0<c<C$, is a solution of the system
(\ref{v2})-(\ref{v0}), then such solution is unique and
coincides with the triple $(V(0),V(1),V(2)).$
\ee{Th}
{\it Proof.}  Let $Y(2)$ be a bounded strictly positive solution of
(\ref{v2}) and let
$$
\int_0^t\psi_u(2)d\widehat M_u+L_t(2)
$$
be the martingale part of $Y(2)$.

Since $Y(2)$ solves
 (\ref{v2}), it follows from the It\^o
formula that for any $\pi\in\Pi(G)$ the process
\begin{equation}\label{ub}
Y^\pi_t=Y_t(2)(1+\int_s^t\pi_ud\widehat S_u)^2
+\int_s^t\pi_u^2(1-\rho^2_u)d\langle M\rangle_u,
\end{equation}
$t\ge s,$ is a local submartingale.

Since $\pi\in\Pi(G)$,
from Lemma \ref{l2.1}
and the Doob inequality we have
$$
E\sup_{t\le T}(1+\int_0^t\pi_ud\widehat S)^2\le
$$
\begin{equation}\label{d}
\le const\big(1+E\int_0^T\pi_u^2\rho_u^2d\langle M\rangle_u
+E\big(\int_0^T|\pi_u\lambda_u|d\langle M\rangle_u\big)^2<\infty
\end{equation}

Therefore, taking in mind that $Y(2)$ is bounded and $\pi\in\Pi(G)$
we obtain that
$$
E\big(\sup_{s\le u\le T}Y_u^\pi\big)^2<\infty
$$
which implies that $Y^\pi\in D$. Thus $Y^\pi$
is a submartingale (as a local submartingale from the class $D$) and by the
boundary condition $Y_T(2)=1$ we obtain
$$
 Y_s(2)\le E\big((1+\int_s^T\pi_ud\widehat S_u)^2+
 \int_s^T\pi^2_u(1-\rho^2_u)d\langle M\rangle_u|G_s\big)
$$
for all $\pi\in\Pi(G)$
and hence
\begin{equation}\label{lev}
Y_t(2)\le V_t(2).
\end{equation}
Let
$$
\tilde\pi_t=-\frac{\lambda_tY_t(2)+\psi_t(2)\rho^2_t}{1-\rho^2_t+\rho^2_t
Y_t(2)}
{\mathcal E}_{t}\left(-\frac{\lambda Y(2)+\psi(2)\rho^2}{1-\rho^2
+\rho^2Y(2)}\cdot\widehat S\right).
$$
Since   $1+\int_0^t\tilde\pi_ud\widehat S_u=
{\mathcal E}_t(-\frac{\lambda Y(2)+\psi(2)\rho^2}
{1-\rho^2+\rho^2Y(2)}\cdot\widehat S)$, it follows
from (\ref{v2}) and the It\^o formula that the
process
$Y^{\tilde\pi}$ defined by (\ref{ub}) is a positive local martingale
and
hence a supermartingale. Therefore
\begin{equation}\label{ley}
Y_s(2)\ge
 E\big((1+\int_s^T\tilde\pi_ud\widehat S_u)^2+
 \int_s^T\tilde\pi^2_u(1-\rho^2_u)d\langle M\rangle_u|G_s\big).
\end{equation}
Let us show that $\tilde\pi$ belongs to the class $\Pi(G)$.

From (\ref{ley}) and (\ref{lev}) we have
for every $s\in [0,T]$
\begin{equation}\label{leo}
 E\big((1+\int_s^T\tilde\pi_ud\widehat S_u)^2+
 \int_s^T\tilde\pi^2_u(1-\rho^2_u)d\langle M\rangle_u|G_s\big)
\le Y_s(2)\le V_s(2)\le1
\end{equation}
and hence
 \begin{equation}\label{leov}
 E\big(1+\int_0^T\tilde\pi_ud\widehat S_u)^2\le1,
\end{equation}
\begin{equation}\label{leop}
E\int_0^T\tilde\pi^2_u(1-\rho^2_u)d\langle M\rangle_u\le1.
\end{equation}
By D*)  the minimal martingale measure $\wh Q^{min}$ for $\wh S$ satisfies the Reverse-H\"older
condition and hence all
conditions of Proposition \ref{p2.1} are satisfied .
Therefore the norm
$$
E\big(\int_0^T\td\pi^2_s\rho_s^2d\langle M\rangle_s\big)+
E\big(\int_0^T|\td\pi_s\ld_s|d\langle M\rangle_s\big)^2
$$
is estimated by $E\big(1+\int_0^T\tilde\pi_ud\widehat S_u)^2$
and hence
$$
E\int_0^T\tilde\pi_u^2\rho_u^2d\langle M\rangle_u<\infty,\;
\;\;
E\big(\int_0^T|\td\pi_s\ld_s|d\langle M\rangle_s\big)^2<\infty.
$$
It follows from (\ref{leop}) and the latter inequality that
$\tilde\pi\in\Pi(G)$ and
from (\ref{ley}) we obtain that
$$
Y_t(2)\ge V_t(2),
$$
which together with (\ref{lev}) gives the equality $Y_t(2)=V_t(2)$.

Thus $V(2)$ is a unique bounded strictly positive solution of
equation (\ref{v2}). Besides
\begin{equation}\label{psi}
\int_0^t\psi(2)_ud\widehat M_u= \int_0^t\varphi(2)_ud\widehat M_u,
\;\;\; L_t(2)=m_t(2)
\end{equation}
for all $t$, $P$-a.s.

Let $Y(1)$ be a solution of equation (\ref{v1}) such that $Y^2(1)\in D$.
By the It\^o formula the process
$$
R_t=Y_t(1){\cal E}_{t}(-\frac{\varphi(2)\rho^2+\lambda V(2)}
{1-\rho^2+\rho^2V(2)}\cdot\widehat S)+
$$
\begin{equation}
\int_0^t{\cal E}_{u}(-\frac{\varphi(2)\rho^2+\lambda V(2)}
{1-\rho^2+\rho^2V(2)}
\cdot\widehat S)\frac{\tilde h_u}{1-\rho_u^2+\rho_u^2V_u(2)}d\langle M\rangle_u
\end{equation}
is a local martingale. Let us show that $R_t$ is a martingale.

 As it was already shown the strategy
$$
\tilde\pi=\frac{\psi_u(2)\rho_u^2+\lambda_u Y_u(2)}
{1-\rho^2+\rho^2Y_u(2)}{\cal E}_{t}(-\frac{\psi(2)\rho^2+\lambda Y(2)}
{1-\rho^2+\rho^2Y(2)}\cdot\widehat S)
$$
belongs to the
class $\Pi(G)$.

Therefore, (see (\ref{d}))
\begin{equation}\label{supe}
E\sup_{t\le T}{\cal E}^2_{t}(-\frac{\psi(2)\rho^2+\lambda Y(2)}
{1-\rho^2+\rho^2Y(2)}\cdot\widehat S)=
E\sup_{t\le T}(1+\int_0^t\tilde\pi_ud\widehat S)^2<\infty
\end{equation}
and hence
$$
Y_t(1){\cal E}_{t}(-\frac{\varphi(2)\rho^2+\lambda V(2)}
{1-\rho^2+\rho^2V(2)}\cdot\widehat S)\in D.
$$
On the other hand equation (\ref{supe}), Lemma \ref{l3.1}  and Lemma \ref{l3.2} imply that
$$
E\sup_{t\le T}\int_0^t{\cal E}_{u}(-\frac{\varphi(2)\rho^2+\lambda V(2)}
{1-\rho^2+\rho^2V(2)}
\cdot\widehat S)\frac{\tilde h_u}{1-\rho_u^2+\rho_u^2V_u(2)}d\langle M\rangle_u\le
$$
$$
\le \frac{1}{c}E\int_0^T{\cal E}_{u}(-\frac{\psi(2)\rho^2+\lambda Y(2)}
{1-\rho^2+\rho^2Y(2)}\cdot\widehat S)
|\tilde h_u|d\langle M\rangle_u
$$
$$
\le\frac{1}{c}E^{1/2}\sup_{t\le T}{\cal E}^2_{t}(-\frac{\psi(2)\rho^2+\lambda Y(2)}
{1-\rho^2+\rho^2Y(2)}\cdot\widehat S)
E^{1/2}\int_0^T\tilde h^2_ud\langle M\rangle_u<\infty.
$$

 Therefore, the process $R_t$ belongs to the class $D$ and hence it
 is a true martingale.
 Using the martingale property and the boundary condition we obtain that
 $$
Y_t(1)=
E\big(\wh{H}_T{\cal E}_{tT}(-\frac{\varphi(2)\rho^2+\lambda V(2)}{1-\rho^2+\rho^2V(2)}
\cdot\widehat S)+
$$
\begin{equation}\label{lin}
+\int_t^T{\cal E}_{tu}(-\frac{\varphi(2)\rho^2+\lambda V(2)}{1-\rho^2+\rho^2V(2)}
\cdot\widehat S)\frac{\tilde h_u}{1-\rho_u^2+\rho_u^2V_u(2)}d\langle M\rangle_u|G_t\big).
\end{equation}
Thus, any solution of (\ref{v1}) is expressed explicitly in terms of $(V(2),\varphi(2))$
in the form  (\ref{lin}). Hence the solution of  (\ref{v1}) is unique and
it coincides with $V_t(1)$.

It is evident that the solution of (\ref{v0}) is also unique.\qed

%{\bf Corollary 4.1}
\begin{cor}
 In addition to conditions A)-C) assume that $\rho$ is a constant
and the mean-variance tradeoff $\la\lambda\cdot M\ra_T$ is deterministic.
Then the solution of (\ref{v2}) is the triple
$(Y(2),\psi(2), L(2))$, with $\psi(2)=0, L(2)=0$ and
\begin{equation}\label{nnu}
Y_t(2)=V_t(2)=\nu(\rho, 1-\rho^2+\la\lambda\cdot M\ra_{T}
-\la\lambda\cdot M\ra_{t}),
\end{equation}
where $\nu(\rho,\alpha)$ is the root of the equation
\begin{equation}\label{fx}
\frac{1-\rho^2}{x}-\rho^2\ln x=\alpha.
\end{equation}
Besides
$$
Y_t(1)=
E\big(H{\cal E}_{tT}(-\frac{\lambda V(2)}{1-\rho^2+\rho^2V(2)}
\cdot\widehat S)
$$
\begin{equation}\label{linn}
+\int_t^T{\cal E}_{tu}(-\frac{\lambda V(2)}{1-\rho^2+\rho^2V(2)}
\cdot\widehat S)\frac{\tilde h_u}{1-\rho^2+\rho^2V_u(2)}d\langle M\rangle_u|G_t\big).
\end{equation}
uniquely solves equation (\ref{v1}) and
the  optimal  filtered wealth process
satisfies the
linear equation
$$
\widehat X_t^*=x-\int_0^t\frac{\lambda_uV_u(2)}
{1-\rho^2+\rho^2V_u(2)}\widehat X^*_ud\widehat S_u
$$
\begin{equation}\label{capi}
+\int_0^t \frac{\varphi_u(1)\rho^2+\lambda_uV_u(1)-\tilde h_u}
{1-\rho^2+\rho^2V_u(2)}
 d\widehat S_u.
\end{equation}
\end{cor}
{\it Proof.} The function $f(x)=\frac{1-\rho^2}{x}-\rho^2\ln x$ is differentiable, strictly
decreasing on $]0,\infty[$ and takes all values from $]-\infty,+\infty[$. So
equation (\ref{fx}) admits a unique solution for all $\alpha$. Besides the inverse function
$\alpha(x)$ is differentiable. Therefore $Y_t(2)$ is a process  of finite variation and it is
adapted since $\la\ld\cdot M\ra_T$ is deterministic.

By definition of $Y_t(2)$ we have that for all $t\in[0,T]$
$$
\frac{1-\rho^2}{Y_t(2)}-\rho^2\ln Y_t(2)=1-\rho^2+
\la\lambda\cdot M\ra_{T}
-\la\lambda\cdot M\ra_{t}.
$$
It is evident that for $\alpha=1-\rho^2$ the solution of
(\ref{fx}) is equal to $1$ and it follows from (\ref{nnu}) that
$Y(2)$ satisfies the boundary conditione $Y_T(2)=1$. Therefore
$$
\frac{1-\rho^2}{Y_t(2)}-\rho^2\ln Y_t(2)-(1-\rho^2)
$$
$$
=-(1-\rho^2)\int_t^Td\frac{1}{Y_u(2)}+\rho^2\int_t^Td\ln Y_u(2)
$$
$$
=\int_t^T\left(\frac{1-\rho^2}{Y^2_u(2)}+\frac{\rho^2}{Y_u(2)}\right)dY_u(2)
$$
and
$$
\int_t^T\frac{1-\rho^2+\rho^2Y_u(2)}{Y_u^2(2)}dY_u(2)=
\la\lambda\cdot M\ra_{T}
-\la\lambda\cdot M\ra_{t}
$$
for all $t\in[0,T]$. Hence
$$
\int_0^t\frac{1-\rho^2+\rho^2Y_u(2)}{Y_u^2(2)}dY_u(2)=
\la\lambda\cdot M\ra_{t}
$$
and integrating both parts of this equality with respect to
 $Y(2)/(1-\rho^2+\rho^2Y(2))$ we obtain that $Y(2)$ satisfies
equation
\begin{equation}\lbl{rrr}
Y_t(2)=Y_0(2)+\int_0^t\frac{Y_u^2(2)\lambda^2_u}{1-\rho^2+
\rho^2Y_u(2)}d\la M\ra_u,
\end{equation}
which implies that the triple $(Y(2),\psi(2)=0, L(2)=0$) satisfies
equation (\ref{v2}) and $Y(2)=V(2)$ by Theorem 4.2.
Equations (\ref{linn}) and (\ref{capi}) follow from (\ref{lin}) and (\ref{cap})
respectively, taking $\varphi(2)=0$.\qed

{\bf Remark 4.1.} In the case of complete information, $M=\widehat M$ and $\rho=1$. Therefore
 equation
(\ref{rrr}) is linear and $Y(2)=e^{\la\lambda\cdot M\ra_{t}
-\la\lambda\cdot M\ra_{T}}$.

{\bf Remark 4.2.} Finally let us make a comment on condition B). It would be desirable to replace condition B) by requiring
that any $G$-martingale is a $F$-semimartingale, but up to now we can't do this, although
 one can weaken this condition imposing that any $G$-martingale is a $\sigma(F^S\vee G)$- martingale, where $\sigma(F_t^S\vee G_t)$
is the minimal $\sigma-$algebra containing $F^S$ and $G_t$, which is satisfied if $F_t^S\subseteq G_t$.

\

\section{Diffusion market model}

\

Let us consider  the financial market model
$$
d\td S_t=\td S_t\mu_t(Y)dt+\td S_t\sigma_t(Y)dw_t^0,
$$
$$
dY_t=a_t(Y)dt+b_t(Y)dw_t,
$$
subjected to initial conditions,
where only the second component $Y$ is observed.
Here $w^0$ and $w$ are corelated Brownian motions with
$Edw^0_tdw_t=\rho dt, \rho\in(-1,1)$.

Let us write
$$
w_t=\rho w_t^0+\sqrt{1-\rho^2}w_t^1,
$$
where  $w^0$ and $w^1$ are independent Brownian motions.
It is evident that $w^\bot=-\sqrt{1-\rho^2}w^0+\rho w^1$ is a
Brownian motion independent of $w$  and one can express Brownian
motions $w^0, w^1$ in terms of $w$ and $w^\bot$ as
\begin{equation}\label{bm}
w_t^0=\rho w_t-\sqrt{1-\rho^2}w_t^\bot, \;\;
w_t^1=\sqrt{1-\rho^2}w_t+\rho w_t^\bot.
\end{equation}

We  assume that $b^2>0,\;\sigma^2>0$ and coefficients $\mu, \sigma, a$ and $b$ are such that
$F_t^{S,Y}=F_t^{w^0,w}$,  $F^Y_t=F_t^w$.
So the stochastic basis will be $(\Omega,{\cal F}, F_t, P)$, where $F_t$
is the natural filtration of $(w^0, w)$ and the flow of observable
events is $G_t=F^w_t$.

Also denote $dS_t=\mu_tdt+\sigma_tdw^0_t$, so that
$d\tilde S_t=\tilde S_tdS_t$
and $S$ is the return of the stock.

Let $\td\pi_t$ be the number shares of the stock
at time $t$. Then $\pi_t=\td\pi_t\td S_t$ represents an amount of money
invested in the stock at the time $t\in[0,T]$.
We consider the mean variance hedging problem
\begin{equation}\lbl{Imv}
\text{to minimize}\;\;\; E[ (x+\int_0^T\td\pi_td\td S_t-H)^2]\;\;\;\;\text{over
all}\;\;\td\pi\;\;\text{for which}\;\;\td\pi\td S\in\Pi(G),
\end{equation}
which is equivalent to study the mean variance hedging problem
$$
\text{to minimize}\;\;\;\;\; E[ (x+\int_0^T\pi_td S_t-H)^2]\;\;\;\;\text{over
all}\;\;\;\;\pi\in\Pi(G).
$$
{\bf Remark 5.1.} Since $S$ is not $G-$adapted,
$\wt\pi_t$ and $\wt\pi_t\wt S_t$ can not be simultaneously  $G$-predictable
 and the problem
\begin{equation}\lbl{Imv1}
\text{to minimize}\;\;\;\;\; E[ (x+\int_0^T\td\pi_td\td S_t-H)^2]\;\;\;\;\text{over
all}\;\;\;\;\td\pi\in\Pi(G),
\end{equation}
is not equivalent to the problem ({\ref{Imv}) and it needs separate consideration.

Comparing with (\ref{str}) we get that in this case
$$
M_t=\int_0^t\sigma_sdw_s^0,\;\;\;\langle M\rangle_t=\int_0^t\sigma_s^2ds,\;\;\;
\ld_t=\frac{\mu_t}{\sigma_t^2}.
$$
It is evident that $w$ is a Brownian motion also with respect to the
filtration $F^{w^0,w^1}$ and condition B) is satisfied. Therefore
by Proposition \ref{p2.2}
$$
\widehat
M_t=\rho\int_0^t\sigma_sdw_s.
$$
By the integral representation theorem the GKW decompositions
(\ref{hhh2}), (\ref{htg}) take following forms
\begin{equation}\label{hhh3}
c_H=EH,\;\;H_t=c_H+\int_0^t h_s\sg_sdw_s^0+\int_0^th_s^1
dw_s^1,
\end{equation}
\begin{equation}\label{htg3}
H_t=c_H+\rho\int_0^t h_s^{G}\sigma_sdw_s+\int_0^th_s^\bot
dw_s^\bot.
\end{equation}
Putting expressions (\ref{bm}) for $w^0, w^1$ in (\ref{hhh3})
and equalizing integrands of (\ref{hhh3}) and (\ref{htg3}) we obtain that
$$
h_t=\rho^2h^G_t-\sqrt{1-\rho^2}\frac{h_t^\bot}{\sigma_t}
$$
and hence
$$
\widehat h_t=\rho^2\widehat{h_t^G}-\sqrt{1-\rho^2}\;
\frac{\widehat{h}_t^\bot}{\sigma_t}.
$$
Therefore by definition of $\widetilde h$
\begin{equation}\lbl{h}
\widetilde h_t=\rho^2
\widehat{h_t^G}-\widehat h_t
=\sqrt{1-\rho^2}\;\frac{\widehat{h}_t^\bot}{\sg_t}
\end{equation}
Using notations
$$
Z_s(0)=\rho\sg_s\vp_s(0),\; Z_s(1)=
\rho\sg_s\vp_s(1),\;Z_s(2)=\rho\sg_s\vp_s(2),\;\;
\theta_s=\frac{\mu_s}{\sg_s}
$$
we obtain the following corollary of Theorem \ref{t4.1}
\be{cor}
Let $H$ be a square integrable $F_T$-measurable
random variable. Then the processes $V_t(0), V_t(1)$ and $V_t(2)$
from
(\ref{ut44}) satisfy the
following system of backward equations
\begin{equation}\label{vw2}
V_t(2)= V_0(2)+\int_0^t \frac{\big(\rho Z_s(2)+
\theta_s V_s(2)\big)^2}
{1-\rho^2+\rho^2V_s(2)}ds
+\int_0^t Z_s(2)dw_s\;\;\;\;V_T(2)=1,
\end{equation}
$$
V_t(1)= V_0(1)+\int_0^t \frac{\big(
\rho Z_s(2)+\theta_sV_s(2)\big)
\big(\rho Z_s(1)+\theta_sV_s(1)-\sqrt{1-\rho^2}\;\widehat{h}_s^\bot\big)}
{1-\rho^2+\rho^2V_s(2)}ds
$$
\begin{equation}\label{vw1}
+\int_0^t Z_s(1)dw_s,\;\;\;\;V_T(1)=E(H|G_T),
\end{equation}
$$
V_t(0)= V_0(0)+\int_0^t \frac{
\big(\rho Z_s(1)+\theta_sV_s(1)-\sqrt{1-\rho^2}\;\widehat{h}_s^\bot\big)^2}
{1-\rho^2+\rho^2V_s(2)}ds
$$
\begin{equation}\label{vw0}
+\int_0^t Z_s(0)dw_s,\;\;V_T(0)=E^2(H|G_T),
\end{equation}
Besides the optimal wealth process $\widehat X^{*}$ satisfies the
linear equation
$$
\widehat X_t^*=x-\int_0^t\frac{\rho Z_s(2)+\theta_sV_s(2)}
{1-\rho^2+\rho^2V_s(2)}\widehat X^*_s(\theta_sds+\rho dw_s)
$$
\begin{equation}\label{capw}
+\int_0^t\frac{\rho Z_s(1)+\theta_sV_s(1)-\sqrt{1-\rho^2}\;\widehat{h}_s^\bot}
{1-\rho^2+\rho^2V_s(2)}(\theta_sds+\rho dw_s).
\end{equation}
\ee{cor}

{\bf Example}. Suppose that $\theta_t$ and $\sg_t$ are deterministic.
Then the solution
of (\ref{vw2}) is the pair
$(V_t(2),\varphi_t(2))$, where $\varphi(2)=0$ and $V(2)$ satisfies
the ordinary differential equation
\begin{equation}\lbl{ode}\frac{dV_t(2)}{dt}=\frac{\theta_t^2V_t^2(2)}
{1-\rho^2+\rho^2V_t(2)},\;\;\;V_T(2)=1.\end{equation}
Solving this equation we obtain that
\begin{equation}\lbl{nu}
V_t(2)=\nu(\rho,1-\rho^2+\int_t^T\theta_s^2ds)\equiv \nu_t^{\theta,\rho},
\end{equation}
where $\nu(\rho,\alpha)$ is the solution of (\ref{fx}).
From (\ref{ode}) follows that
\begin{equation}\lbl{mart}
(\ln\nu_t^{\theta,\rho})'=\frac{\theta_t^2\nu_t^{\theta,\rho}}{1-\rho^2+\rho^2\nu_t^{\theta,\rho}}\;\text{and}\;
\ln\frac{\nu_s^{\theta,\rho}}{\nu_t^{\theta,\rho}}
=\int_t^s\frac{\theta_r^2\nu_r^{\theta,\rho}dr}{1-\rho^2+\rho^2\nu_r^{\theta,\rho}}.
\end{equation}
If we solve the linear BSDE (\ref{vw1}) and use (\ref{mart}) we obtain
\begin{align*}
V_t(1)=E\left[\widehat{H}_T(w){\cal E}_{tT}\left(-\int_0^\cdot\frac{\theta_r\nu_r^{\theta,\rho}}{1-\rho^2+\rho^2\nu_r^{\theta,\rho}}(\theta_rdr+\rho dw_r)\right)|{G}_t\right]\\
\int_t^T\frac{\theta_s\nu_s^{\theta,\rho}\sg_s}{1-\rho^2+\rho^2\nu_s^{\theta,\rho}}E\left[\td{h}_s(w){\cal E}_{ts}\left(-\int_0^\cdot\frac{\theta_r\nu_r^{\theta,\rho}}{1-\rho^2+\rho^2\nu_r^{\theta,\rho}}(\theta_rdr+\rho dw_r)\right)|{G}_t\right]ds\\
=\nu_t^{\theta,\rho}E\left[\widehat{H}_T(w){\cal E}_{tT}\left(-\int_0^\cdot\frac{\theta_r\nu_r^{\theta,\rho}}{1-\rho^2+\rho^2\nu_r^{\theta,\rho}}\rho dw_r\right)|
{G}_t\right]\\
+\nu_t^{\theta,\rho}\int_t^T\frac{\mu_s}{1-\rho^2+\rho^2\nu_s^{\theta,\rho}}E\left[\td{h}_s(w){\cal E}_{ts}\left(-\int_0^\cdot\frac{\theta_r\nu_r^{\theta,\rho}}{1-\rho^2+\rho^2\nu_r^{\theta,\rho}}\rho dw_r\right)|{G}_t\right]ds
\end{align*}
Using the Girsanov theorem we  finally get
\begin{align}\notag
V_t(1)=\nu_t^{\theta,\rho}
E\left[\widehat{H}_T\left(\int_t^\cdot\rho\frac{\theta_r\nu_r^{\theta,\rho}}{1-\rho^2+\rho^2\nu_r^{\theta,\rho}} dr+w\right)|{G}_t\right]\\
\label{vv1}
+\nu_t^{\theta,\rho}\int_t^T\frac{\mu_s}{1-\rho^2+\rho^2\nu_s^{\theta,\rho}}
E\left[\td{h}_s\left(\rho\int_t^\cdot\frac{\theta_r\nu_r^{\theta,\rho}}{1-\rho^2+\rho^2\nu_r^{\theta,\rho}} dr+w\right)|
{G}_t\right]ds.
\end{align}
\qed

Suppose now that $H=H(w_T^0,w_T),\;Y=w$ and
$\frac{\mu_t}{\sigma_t}=\theta(t,w_t)$ for some continuous
function $\theta(t,x)$ and a differentiable function $H(x,y)$.
Then using the elementary ideas of Malliavin's calculus we get
$$\rho h_t^{G}\sg_t\equiv h^{\cal G}(t,w_t^0,w_t)=E[\rho \partial_xH(w_T^0,w_T)+\partial_yH(w_T^0,w_T)|{\cal F}_t],$$
$$
h_t^\bot=h^\bot(t,w_t^0,w_t)=-\sqrt{1-\rho^2}E[\partial_xH(w_T^0,w_T)|{\cal F}_t],$$
where $\partial_x{H}(x,y),\partial_y{H}(x,y)$ denote the partial derivatives of $\cal H$.
It is evident that \\
$E[f(t,w_t^0,w_t)|{\cal F}_t^w]=
E[f(t,\rho w_t-\sqrt{1-\rho^2}w_t^\bot, w_t)|w_t]$ for $f=h^w,h^\bot,H$.
Thus we obtain the exact expression for $\widehat{H}_T(y),
{\widehat h}^{\cal G}(t,y)$ and $\widehat{h}^\bot(t,y)$
\begin{align}\lbl{mal}
\widehat{H}_T(y)=EH(\rho y-\sqrt{1-\rho^2}w_T^\bot,y)\equiv E[H(w_T^0,y)|w_T=y],\\
\notag
{\widehat h}^{\cal G}(t,y)=Eh^{\cal G}(t,\rho y-\sqrt{1-\rho^2}w_t^\bot, y),\\
\notag
{\widehat h}^\bot(t,y)=Eh^\bot(t,\rho y-\sqrt{1-\rho^2}w_t^\bot, y)=-\sqrt{1-\rho^2}E[H_x(w_T^0,w_T)|w_t=y]\\
\notag
\equiv-\sqrt{1-\rho^2}E\partial_xH(\rho (y+w_T-w_t)-\sqrt{1-\rho^2}w_T^\bot,\rho(y+w_T-w_t)).
\end{align}

{\bf Remark 5.2.} For deterministic $\theta_t$ the equalities
\begin{align*}E\left[\widehat{H}_T
\left(\rho\int_t^T\frac{\theta_r\nu_r^{\theta,\rho}}{1-\rho^2+\rho^2\nu_r^{\theta,\rho}} dr+w_T\right)|w_t=y\right]
\\
=E\bigg[{H}\bigg(\rho\left(\rho\int_t^T\frac{\theta_r\nu_r^{\theta,\rho}dr}{1-\rho^2+\rho^2\nu_r^{\theta,\rho}} +w_T\right)-\sqrt{1-\rho^2}\;w_T^\bot,\\
\rho\int_t^T\frac{\theta_r\nu_r^{\theta,\rho}dr}{1-\rho^2+\rho^2\nu_r^{\theta,\rho}}+w_T\bigg)|w_t=y\bigg]
\\
\equiv E{H}\bigg(\rho\left(\rho\int_t^T\frac{\theta_r\nu_r^{\theta,\rho}dr}{1-\rho^2+\rho^2\nu_r^{\theta,\rho}} +w_T-w_t+y\right)-\sqrt{1-\rho^2}\;w_T^\bot,\\
\rho\int_t^T\frac{\theta_r\nu_r^{\theta,\rho}dr}{1-\rho^2+\rho^2\nu_r^{\theta,\rho}}+w_T-w_t+y\bigg)
\end{align*}
and
\begin{align*}E\left[{\wh h}_s^\bot\left(\rho\int_t^s\frac{\theta_r\nu_r^{\theta,\rho}}{1-\rho^2+\rho^2\nu_r^{\theta,\rho}} dr+w_s\right)|w_t=y\right]
\\
=-\sqrt{1-\rho^2}E\bigg[E\bigg[{\partial_xH}\bigg(\rho\left(\rho\int_t^s\frac{\theta_r\nu_r^{\theta,\rho}dr}
{1-\rho^2+\rho^2\nu_r^{\theta,\rho}} +w_T\right)-\sqrt{1-\rho^2}\;w_T^\bot,\\
\rho\int_t^s\frac{\theta_r\nu_r^{\theta,\rho}dr}{1-\rho^2+\rho^2\nu_r^{\theta,\rho}}+w_T\bigg)|w_s\bigg]|w_t=y\bigg]
\\
=-\sqrt{1-\rho^2}E\bigg[{\partial_xH}\bigg(\rho\left(\rho\int_t^s\frac{\theta_r\nu_r^{\theta,\rho}dr}
{1-\rho^2+\rho^2\nu_r^{\theta,\rho}} +w_T\right)-\sqrt{1-\rho^2}\;w_T^\bot,\\
\rho\int_t^s\frac{\theta_r\nu_r^{\theta,\rho}dr}{1-\rho^2+\rho^2\nu_r^{\theta,\rho}}+w_T\bigg)|w_t=y\bigg]
\\
\equiv -\sqrt{1-\rho^2}E{\partial_xH}\bigg(\rho\left(\rho\int_t^T
\frac{\theta_r\nu_r^{\theta,\rho}dr}{1-\rho^2+\rho^2\nu_r^{\theta,\rho}} +w_T-w_t+y\right)-\sqrt{1-\rho^2}\;w_T^\bot,\\
\rho\int_t^T\frac{\theta_r\nu_r^{\theta,\rho}dr}{1-\rho^2+\rho^2\nu_r^{\theta,\rho}}+w_T-w_t+y\bigg)
\end{align*}
are valid.
\qed

Using the well known connection between BSDEs and PDEs  we can prove
the following

\be{prop}\label{pr3}  The system of nonlinear PDEs
\begin{gather}\label{pde}
\partial_t v(2)+\frac{1}{2}\partial_{y}^2v(2)=\frac{(\theta(t,y)v(2)+\rho\partial_y v({2}))^2}{1-\rho^2+\rho^2 v(2)},\;\;\;\;v(2,T,y)=1,\\\label{pde2}
\partial_t v(1)+\frac{1}{2}\partial_{y}^2v(1)\\
\notag
=\frac{\left((\theta(t,y)v(2)+\rho\partial_y v({2})\right)\left((\theta(t,y)v(1)+\rho\partial_y v(1))
+(1-\rho^2)E[\partial_xH(w_T^0,w_T)|w_t=y]\right)}
{1-\rho^2+\rho^2 v(2)},\\
\notag v(1,T,y)=E[H(w_T^0,y)|w_T=y]
\end{gather}
admits sufficiently smooth solution  and the solution of (\ref{v2}),(\ref{v1}) can be
represented as $V_t(1)=v(1,t,w_t), Z_{t}(1)=\partial_y v(1,t,w_t), V_t(2)=v(2,t,w_t),Z_{t}(2)=
\partial_{y}v(2,t,w_t)$.

Besides the optimal strategy is of the form
\begin{align}\label{stra}
\pi^*(t,X,y)=-\frac{\theta(t,y)v(2,t,y)+\rho
\partial_y v(2,t,y)}{1-\rho^2+\rho^2 v(2,t,y)}X\\
\notag
+\frac{\left(\theta(t,y)v(1,t,y)+\rho \partial_y v(1,t,y)
+(1-\rho^2)E[\partial_xH(w_T^0,w_T)|w_t=y]\right)}
{1-\rho^2+\rho^2 v(2,t,y)},
\end{align}
where $X$ and $y$ are states at time $t$ of the wealth and of an
observable process.
\ee{prop}

{\bf Example(continued)}. We suppose in addition that $\theta,\sg$
 are constants
and $H={\cal H}(S_T,Y_T)\equiv {\cal H}(\mu T+\sg w_T^0,Y_T)$.  Then  using the equality $\frac{1}{\theta}\ln\frac{\nu_s^{\theta,\rho}}{\nu_t^{\theta,\rho}}=
\int_t^s\frac{\theta\nu_r^{\theta,\rho}dr}{1-\rho^2+\rho^2\nu_r^{\theta,\rho}}$
the formula (\ref{vv1}) can be simplified.
It is easy to see that $$\lim_{\theta\to 0}\frac{1}{\theta}\ln\nu_t^{\theta,\rho}=0.$$
Thus we can set that expression $\frac{1}{\theta}\ln\nu_t^{\theta,\rho}$  is zero as $\theta=0$.
For simplicity we also assume
that $a=0,\;b=1$. For $v(1)$ using (\ref{vv1}) and Remark 5.2 we get
\begin{align*}
v(1,t,y)=\nu_t^{\theta,\rho}E\left[\widehat{H}\left(-\frac{\rho}{\theta}\ln\nu_t^{\theta,\rho}
+w_T-w_t+y\right)\right]\\
+\nu_t^{\theta,\rho}\int_t^T\frac{\mu}{1-\rho^2+\rho^2\nu_s^{\theta,\rho}}E\left[
\td{h}_s\left(\frac{\rho}{\theta}\ln\frac{\nu_s^{\theta,\rho}}{\nu_t^{\theta,\rho}}+w_s-w_t+y\right)\right]ds,
\end{align*}
or equivalently
\begin{gather}\lbl{vm1}
\quad v(1,t,y)=\nu_t^{\theta,\rho}\\
\notag
\times E{\cal H}\left(\mu T+\sg\rho(-\frac{\rho}{\theta}\ln\nu_t^{\theta,\rho}+w_T-w_t+y)-\sg\sqrt{1-\rho^2}w_T^\bot,
\rho(-\frac{\rho}{\theta}\ln\nu_t^{\theta,\rho}+w_T-w_t+y)\right)\\
\notag+(1-\rho^2)\mu\nu_t^{\theta,\rho}\int_t^T\frac{1}{1-\rho^2+\rho^2\nu_s^{\theta,\rho}}\\
\notag\times
E\partial_x{\cal H}\left(\mu T+\sg\rho(\frac{\rho}{\theta}\ln\frac{\nu_s^{\theta,\rho}}{\nu_t^{\theta,\rho}}+w_T-w_t+y)-\sg\sqrt{1-\rho^2}w_T^\bot,
\rho(\frac{\rho}{\theta}\ln\frac{\nu_s^{\theta,\rho}}{\nu_t^{\theta,\rho}}+w_T-w_t+y)\right)ds.
\end{gather}
taking in mind (\ref{mal}).
This formula together with (\ref{nu}) and(\ref{stra}) gives
an explicit solution of the problem (\ref{Imv})
for the case of constant coefficients.

\

\appendix
\section{Appendix}

\

For convenience we give the proofs of the following assertions used in the paper.

%Lemma 2.1.
\be{lem}\lbl{l2.1}
Let conditions {\rm A)--C)} be satisfied and
$\widehat M_t=E(M_t|G_t)$. Then $ \la \widehat M\ra $ is
absolutely continuous w.r.t $ \langle M\rangle $ and $\mu^{<M>}$ a.e.
$$
\rho^2_t=\frac{d \la \widehat M\ra _t}{d \langle M\rangle _t}\le 1
$$

\ee{lem}

{\bf Proof.}
By (\ref{mg})
for any bounded ${G}$-predictable process $h$
\begin{align}\lbl{l4}
\notag
E\int_0^th_s^2d \la \widehat M\ra _s=
E\left(\int_0^th_sd\widehat M_s\right)^2
=E\left(E\left(\int_0^th_sdM_s\big|G_t\right)\right)^2\le\\
\le E\left(\int_0^th_sdM_s\right)^2\le E\int_0^th_s^2d \langle M\rangle_s
\end{align}
which implies that $ \la \widehat M\ra $ is absolutely continuous
w.r.t $ \langle M\rangle $, i.e.,
$$
 \la \widehat M\ra _t=\int_0^t\rho^2_sd \langle M\rangle _s
$$
for a ${G}$-predictable process $\rho$.

Moreover (\ref{l4}) implies that the process $\la M\ra - \la\wh M\ra$ is increasing
and hence $\rho^2\le 1$ $\mu^{\langle M\rangle}$ a.e.

% Lemma 3.1
\be{lem}\lbl{l3.1}
 Let $H\in L^2(P,F_T)$ and let conditions $A)-C)$ be
satisfied. Then
\begin{equation}\label{tdh}
E\int_0^T\tilde h^2_ud\langle M\rangle_u<\infty.
\end{equation}
\ee{lem}
{\it Proof.}
It is evident that
$$
E\int_0^T(h_u^G)^2d\la \widehat M\ra_u <\infty,
\;\;\;\;E\int_0^Th_u^2d\langle M\rangle_u <\infty.
$$
Therefore, by definition of $\tilde h$ and Lemma \ref{l2.1},
$$
E\int_0^T\tilde h^2_ud\langle M\rangle_u
$$
$$
\le 2E\int_0^TE^2(h_u|G_u)d\langle M\rangle_u+
2E\int_0^TE^2(h^G_u|G_u)\rho^4_ud\langle M\rangle_u
$$
$$
\le 2E\int_0^T h^2_ud\langle M\rangle_u+
2E\int_0^T (h^G_u)^2\rho_u^2d\langle \widehat M\rangle_u<\infty.
$$
Thus $\tilde h\in\Pi(G)$ by remark 2.3.\qed

%{\bf Lemma A4.}
\be{lem}\lbl{la4}
a) Let $Y=(Y_t, t\in[0,T])$ be a bounded positive submartingale with
the canonical decomposition
$$
Y_t=Y_0+B_t+m_t
$$
where $B$ is a predictable
increasing process and $m$ is a martingale. Then $m\in BMO$.

b) In particular the martingale part of $V(2)$ belongs to $BMO$.
If $H$ is bounded, then martingale parts of $V(0)$ and $V(1)$
also belong to the class $BMO$, i.e., for $i=0,1,2,$
\begin{equation}\label{psi}
E(\int_\tau^T\varphi^2_u(i)\rho^2_ud\langle M\rangle_u|G_\tau)+
E(\la m(i)\ra_T-\la m(i)\ra_\tau|G_\tau)\le C
\end{equation}
for every stopping time $\tau$.
\ee{lem}

{\it Proof.}
Applying the It\^o formula for $Y_T^2-Y^2_\tau$ we have
\begin{equation}\label{cond}
\langle m\rangle_T-\langle m\rangle_\tau+2\int_\tau^TY_udB_u+2\int_\tau^TY_udm_u=Y^2_T-Y^2_\tau\le
 Const.
\end{equation}
Since $Y$ is positive and $B$ is an increasing process, taking conditional
expectations in (\ref{cond}) we obtain
$$
E( \langle m\rangle_T-\langle m\rangle_\tau|F_\tau)\le Const.
$$
for any stopping time $\tau$, hence $m\in BMO$.

(\ref{psi}) follows from assertion a) applied for positive submartingales
$V(0), V(2)$ and $V(0)+V(2)-2V(1)$. For the case $i=1$ one should
take into account also the inequality
$$
\la m(1)\ra_t \le const(\la m(0)+m(2)-2m(1)\ra_t + \la m(0)\ra_t +\la m(2)\ra_t).
$$

{\bf Acknowledgments.}
We would like to thank N. Lazrieva for useful remarks and comments.


\begin{thebibliography}{120}

\bibitem{B} J. M. Bismut, Conjugate convex functions in optimal stochastic
control, {\it J. Math. Anal. Appl.} {\bf 44} (1973), 384--404.

\bibitem{Ch}  R. Chitashvili, Martingale ideology in the theory of
controlled stochastic processes, {\it Lect. Notes in Math.} {\bf
1021}, 73--92, Springer-Verlag, New York, 1983.

\bibitem{DPR} G. B. Di Masi, E. Platen and W. J. Runggaldier , Hedging of Options under
Discrete Observation on Assets with Stochastic Volatility",
Seminar on Stoch. Anal. Rand. Fields Appl.,(1995) 359-364.

\bibitem{D-M-S-S-S}  F. Delbaen, P. Monat, W. Schachermayer, W. Schweizer and C. Stricker,
Weighted norm inequalities and hedging in incomplete
markets, {\it Finance Stoch.} {\bf 1, No. 2}, (1997), 181-227.

\bibitem{D-M}   C. Dellacherie and P. A. Meyer, Probabilit\'{e}s et potentiel, II.
{\em Hermann, Paris,} 1980.

\bibitem{D-R}  D. Duffie and H.R. Richardson, Mean-Variance hedging in continuous time.
{\em Ann. Appl. Probab.} {\bf 1}(1991), 1--15.

\bibitem{El-Q} N.El Karoui and M.C. Quenez, Dynamic programming and
pricing of contingent claims in an incomplete market, {\it SIAM J. Control Optim.}
{\bf 33} 1 (1995), 29-66.

\bibitem{F-S}  H. F\"ollmer and D. Sondermann,
Hedging non-redundant contingent claims. {\em Contributions to
mathematical economics hon. G. Debreu $($W. Hildenbrand and A.
Mas-Collel, eds.$)$,}. {\em North Holland, Amsterdam,} (1986),
205--223.

\bibitem{FR}
R. Frey and W. J. Runggaldier,  Risk-minimizing hedging strategies
under restricted information: the case of stochastic volatility
models observable only at discrete random times. Financial
optimization. Math. Methods Oper. Res. 50 (1999), no. 2, 339--350.


\bibitem{G-L-Ph}  C. Gourieroux, J. P. Laurent, and H. Pham, Mean-variance hedging and
numeraire. {\em  Math. Finance} {\bf 8}, No. 3(1998), 179--200.


\bibitem{H-P-Sc}  D. Heath, E. Platen, and M. Schweizer, A comparision of
two quadratic approaches to hedging in incomplete markets. {\em
Math. Finance } {\bf 11} No. 4 (2001) 385--413.


\bibitem{H}  C. Hipp, Hedging general claims. {\em  Proc. $3$rd AFIR Colloquium}
2, 603--613, {\em Rome,} 1993.

\bibitem{Hof-P-Sc}  N. Hoffman, E. Platen and M. Schweizer,
Option pricing under incompleteness and stochastic volatility,
{\it Math. Finance}, Vol. 2, N. 3, 1992, pp. 153-187.

\bibitem{J}  J. Jacod, Calcule Stochastique et probl\`{e}mes des martingales.
{\em Lecture Notes in Math.} 714, {\em Springer, Berlin etc.,}
(1979).


\bibitem{Kz}  N. Kazamaki, Continuous exponential martingales and {\rm BMO},
{\it Lecture Notes in Math.} {\bf 1579}, Springer, New York, 1994.


\bibitem{Kr} D.O.Kramkov, Optional decomposition of
supermartingales and hedging contingent claims in incomplete
security markets, {\it Probab. Theory Related Fields} {\bf 105}
(1996), 459--479.


\bibitem{Kob} {Kobylanski M.}, (2000) Backward stochastic differential
equation and
partial differential equations with quadratic growth, The Annals of
Probability,
vol. 28, N2, 558-602.

\bibitem{LS} {Lepeltier J.P. and  San Martin J.}, (1998) Existence for
BSDE with
superlinear-quadratic coefficient, Stoch. Stoch. Rep. 63, 227-240.


\bibitem{L-Sh2}  R.Sh. Liptzer and A.N. Shiryayev, Martingale theory,  Nauka,
Moscow, 1986.

\bibitem{MT7} M. Mania and R. Tevzadze, Backward Stochastic PDE and
Imperfect Hedging, International Journal of Theoretical and
Applied Finance, vol.6, 7,(2003),663-692.

\bibitem{Mrl} Morlais M.A., Quadratic Backward Stochastic Differential
Equations (BSDEs) Driven by a Continuous Martingale and
Application to the Utility Maximization Problem,
http://hal.ccsd.cnrs.fr/ccsd-00020254/

\bibitem{Par-P}  E. Pardoux and S.G.Peng, Adapted solution of a backward
stochastic differential equation, {\it Systems Control Lett.} {\bf
14} (1990),55--61.

\bibitem{Ph} H. Pham, Mean-variance hedging for partially observed drift processes
Int. J. Theor. Appl. Finance 4 (2001), no. 2, 263--284.

\bibitem{R-Sc}  T. Rheinl\"ander and M. Schweizer, On $L^2$- Projections
on a space of Stochastic Integrals, {\it Ann. Probab.} {\bf 25}
(1997), 1810--1831.

\bibitem{S} M. Sch\"al, On quadratic cost criteria for option hedging,
{\it Math. Oper. Res.} {\bf 19} (1994), 121--131.

\bibitem{S92}   M. Schweizer, Mean-variance hedging for general claims. {\em Ann. Appl.
Probab.} {\bf 2}(1992), 171--179.

\bibitem{S94}   M. Schweizer, Approximating random variables by stochastic
integrals. {\em Ann. Probab.} {\bf 22},(1994), No. 3, 1536--1575.

\bibitem{Sc94}   M. Schweizer, Risk-minimizing hedging strategies under restricted information.
Math. Finance 4, No.4,  (1994), 327-342.

\bibitem{Tev} R. Tevzadze, Solvability of backward stochastic differential equations with
quadratic growth, to appear in Stochastic Proccesses and their Applications.

\end{thebibliography}
\end{document}